\documentclass[11pt]{article}

\usepackage{amscd,amsmath, amssymb, fancyhdr, epsfig, color, tocloft, mathtools}
\usepackage{dutchcal}

\usepackage[backref=page]{hyperref}
\renewcommand*{\backrefalt}[4]{%
	\ifcase #1 (Not cited.)%
	\or        (Cited on page~#2.)%
	\else      (Cited on pages~#2.)%
	\fi}

\hypersetup{
	colorlinks   = true,
	citecolor    = magenta,
	linkcolor    = blue,
	urlcolor     = magenta	
}

\newcommand{\version}{version 1.0,\ \ Jan. 3, 2022}
\makeatletter
\def\x@arrow{\DOTSB\Relbar}
\def\xlongequalsignfill@{\arrowfill@\x@arrow\Relbar\x@arrow}
\providecommand{\xlongequal}[2][]{%
	\ext@arrow 0099\xlongequalsignfill@{#1}{#2}}
\def\xlongrightarrowfill@{\arrowfill@\relbar\relbar\longrightarrow}
\makeatother

\numberwithin{equation}{section}

\def\eqref#1{(\ref{#1})}
\newcommand{\goth}{\mathfrak}

\newcommand{\Z}{{\mathbb Z}}
\newcommand{\C}{{\mathbb C}}
\newcommand{\R}{{\mathbb R}}

\def\1{\sqrt{-1}\:}
\newcommand{\restrict}[1]{{\left|_{{\phantom{|}\!\!}_{#1}}\right.}}
\newcommand{\cntrct}                
{\hspace{2pt}\raisebox{1pt}{\text{$\lrcorner$}}\hspace{2pt}}
\newcommand{\arrow}{{\:\longrightarrow\:}}

\newcommand{\calo}{{\cal O}}
\newcommand{\cac}{{\cal C}}

\renewcommand{\bar}{\overline}
\renewcommand{\phi}{\varphi}
\renewcommand{\epsilon}{\varepsilon}
\renewcommand{\geq}{\geqslant}

\newcommand{\Tot}{\operatorname{Tot}}
\newcommand{\Id}{\operatorname{Id}}

\newcommand{\Hom}{\operatorname{Hom}}

\newcommand{\Aut}{\operatorname{Aut}}

\newcommand{\codim}{\operatorname{codim}}

\newcommand{\Lie}{\operatorname{Lie}}

\newcommand{\GL}{\operatorname{GL}}

\newcommand{\U}{{\operatorname{U}}}

\newcounter{Mycounter}[section]
\newcounter{lemma}[section]
\renewcommand{\thelemma}{{Lemma \thesection.\arabic{lemma}}}
\newcommand{\lemma}{%
	\setcounter{lemma}{\value{Mycounter}}
	\refstepcounter{lemma}
	\stepcounter{Mycounter}
	{\noindent \bf \thelemma:\ }}

\newcounter{claim}[section]
\newcounter{sublemma}[section]
\newcounter{corollary}[section]
\newcounter{theorem}[section]
\renewcommand{\thetheorem}{{Theorem \thesection.\arabic{theorem}}}
\newcommand{\theorem}{%
	\setcounter{theorem}{\value{Mycounter}}
	\refstepcounter{theorem}
	\stepcounter{Mycounter}
	{\noindent \bf \thetheorem:\ }}

\newcounter{conjecture}[section]
\renewcommand{\theconjecture}{{Conjecture \thesection.\arabic{conjecture}}}
\newcommand{\conjecture}{%
	\setcounter{conjecture}{\value{Mycounter}}
	\refstepcounter{conjecture}
	\stepcounter{Mycounter}
	{\noindent \bf \theconjecture:\ }}

\newcounter{proposition}[section]
\renewcommand{\theproposition} {{Proposition \thesection.\arabic{proposition}}}
\newcommand{\proposition}{%
	\setcounter{proposition}{\value{Mycounter}}
	\refstepcounter{proposition}
	\stepcounter{Mycounter}
	{\noindent \bf \theproposition:\ }}

\newcounter{definition}[section]
\renewcommand{\thedefinition} {{Definition~\thesection.\arabic{definition}}}
\newcommand{\definition}{%
	\setcounter{definition}{\value{Mycounter}}
	\refstepcounter{definition}
	\stepcounter{Mycounter}
	{\noindent \bf \thedefinition:\ }}

\newcounter{example}[section]
\renewcommand{\theexample}{{Example \thesection.\arabic{example}}}
\newcommand{\example}{%
	\setcounter{example}{\value{Mycounter}}
	\refstepcounter{example}
	\stepcounter{Mycounter}
	{\noindent \bf \theexample:\ }}

\newcounter{remark}[section]
\renewcommand{\theremark}{{Remark \thesection.\arabic{remark}}}
\newcommand{\remark}{%
	\setcounter{remark}{\value{Mycounter}}
	\refstepcounter{remark}
	\stepcounter{Mycounter}
	{\noindent \bf \theremark:\ }}

\newcounter{problem}[section]
\newcounter{question}[section]
\makeatletter

\@addtoreset{equation}{section}
\@addtoreset{footnote}{section}
\makeatother

\def\blacksquare{\hbox{\vrule width 5pt height 5pt depth 0pt}}
\def\endproof{\blacksquare}

\newcommand{\proof}{{\bf Proof: \ }}
\newcommand{\pstep}{{\bf Proof. Step 1: \ }}

\begin{document}
	
	\begin{center}
		{\Large\bf  Holomorphic tensors on Vaisman
                manifolds}\\[5mm]
		{\large
			Liviu Ornea\footnote{Liviu Ornea is  partially supported by Romanian Ministry of Education and Research, Program PN-III, Project number PN-III-P4-ID-PCE-2020-0025, Contract  30/04.02.2021},  
			Misha Verbitsky\footnote{Misha Verbitsky is partially supported 
				by the HSE University Basic Research Program, FAPERJ E-26/202.912/2018 
				and CNPq - Process 313608/2017-2.\\[1mm]
				
			}\\[4mm]
			
		}
		
	\end{center}
	
\begin{center}
{\em \small Dedicated to Professor Valentin Poenaru at his ninetieth birthday}
\end{center}
	
	\hfill
	
	{\small
		\hspace{0.15\linewidth}
		\begin{minipage}[t]{0.7\linewidth}
			{\bf Abstract} \\
An LCK (locally conformally K\"ahler)
manifold is a complex manifold admitting a Hermitian form $\omega$
which satisfies $d\omega =\omega\wedge \theta$, where
$\theta$ is a closed 1-form, called {\bf the Lee form}.
An LCK manifold is called {\bf Vaisman} if the Lee
form is parallel with respect to the Levi-Civita connection.
The dual vector field, called {\bf the Lee field}, is
holomorphic and Killing. We prove that any holomorphic
tensor on a Vaisman manifold is invariant with respect
to the Lee field. This is used to compute the Kodaira
dimension of Vaisman manifolds. We prove that the
Kodaira dimension of a Vaisman manifold 
obtained as a $\Z$-quotient of an algebraic cone over a 
projective manifold $X$ is equal to the Kodaira dimension 
of $X$. This can be applied to prove the deformational
stability of the Kodaira dimension of Vaisman manifolds.
		\end{minipage} 
	}

	\tableofcontents
	
	
	\section{Introduction}
	\label{_Intro_Section_}

Locally conformally K\"ahler geometry
is probably the simplest complex geometry beyond 
the K\"ahler which admits many examples, explicit
and non-explicit. A manifold is called {\bf LCK 
(locally conformally K\"ahler)} if it admits
a K\"ahler covering, with the deck group acting
by holomorphic homotheties (\ref{_LCK_Definition_}). 
Every Hopf manifold 
is LCK, because its universal cover admits a K\"ahler
potential, with the deck transform multiplying
the potential by a number. A manifold with
this property is called {\bf an LCK manifold with potential};
an LCK manifold admits a potential if and only
if it can be realized as a complex submanifold
in a Hopf manifold (\cite{ov_lckpot}).

Recall that {\bf a  diagonal
Hopf manifold} is a manifold of
form $\frac{\C^n\backslash 0}{\langle A\rangle}$,
where $A$ is a linear contraction which
can be diagonalized in appropriate basis.
A more special class of LCK manifolds
with potential are {\bf Vaisman manifolds},
Subsection \ref{_Vaisman_Subsubsection_}, which can 
be characterized as complex manifolds
admitting a holomorphic embedding to a diagonal
Hopf manifold (\cite[Theorem 4.9]{_OV:Lee_classes_}).

Every Vaisman manifold $M$ is equipped with an action of
the Lee and anti-Lee fields, which are Killing and holomorphic
(\ref{_Lee_field_Killing_Corollary_}).
In \cite{_Verbitsky:Sta_Elli_}, 
it was shown that any stable holomorphic bundle on $M$
is equivariant with respect to this action. 
From \ref{_Subva_Vaisman_Theorem_}, it 
follows that any complex subvariety is invariant
under this action as well.

Let $G$ be the group of holomorphic isometries
of a Vaisman manifold $M$ obtained as the closure of the
group generated by the Lee and the anti-Lee flow.
In this paper, we prove that any
holomorphic tensor field on $M$ is
$G$-invariant. 

For holomorphic vector fields and differential forms this result was obtained
by K. Tsukada in \cite[Theorem 3.3, Theorem 4.2]{tsuk}.

We prove a similar result for LCK manifolds with potential.
Let $M$ be an LCK manifold with potential, and $\tilde M$
its K\"ahler $\Z$-cover. We assume that $\tilde M$ is an open 
algebraic cone (\ref{_alge_cone_Definition_}); 
this is true automatically when $\dim_\C M\geq 3$
(\cite{ov_ac}).
When $\dim_\C M=2$, this would follow if we assume 
the GSS conjecture 
(\cite[Section 25.6]{_OV:book_}). Since $\tilde M$ is algebraic,
it makes sense to speak of the Zariski closure of the
$\Z$-action on $\tilde M$. We consider the Zariski 
closure
${\cal G}$ of the $\Z$-action as the smallest algebraic
subgroup of $\Aut(\tilde M)$ containing this $\Z$-action;
this group is clearly commutative and positive-dimensional. Since 
the action of ${\cal G}$ commutes with the action of the deck
transform group $\Z\subset {\cal G}$ on $\tilde M$, 
the quotient ${\cal G}/\Z$ naturally acts on $M$.

We prove that any holomorphic tensor field on 
an LCK manifold with potential is ${\cal G}$-invariant
(\ref{_LCK_pot_Zariski_closure_Proposition_}).


Let $H = \frac{\C^n \backslash 0}{\langle A \rangle}$
be a diagonal Hopf manfold, and ${\cal G}$ the
Zariski closure of ${\langle A \rangle}$.
Let $\alpha_i$ be the eigenvalues of $A$,
and $A_1$ an operator which has the eigenvalues $|\alpha_i|$
in the same basis. In  \ref{_Zariski_contains_real_part_Theorem_},
we show that the linear field $\log \!A_1$ belongs to the
Zariski closure of ${\langle A \rangle}$, and 
in \ref{_holo_tensor_on_Vaisman_Lee_invariant_Theorem_},
step 3, we prove that $\log A_1$ can be obtained
as the Lee field of a Vaisman structure on $H$.

This implies that the Lee and anti-Lee field of a Vaisman 
manifold belong to the Lie algebra of ${\cal G}$
(the Zariski closure of the $\Z$-action).
Therefore, any holomorphic tensor field is Lee and anti-Lee
invariant.

\section{Preliminaries}

For the material in this section, we refer the reader to \cite{_OV:book_}.

\subsection{Locally conformally K\"ahler manifolds}

\definition\label{_LCK_Definition_}
A complex  manifold $(M, I)$ is called {\bf locally conformally
	K\"ahler} (LCK, for short) if it admits a
covering $(\tilde M, I)$ equipped with a K\"ahler metric
$\tilde\omega$ such that the deck group of the
cover acts on $(\tilde M, \tilde \omega)$ by holomorphic 
homotheties. An {\bf LCK metric} on an LCK manifold
is an Hermitian metric on $(M,I)$ such that its
pullback to $\tilde M$ is conformal with $\tilde \omega$.

\hfill

\proposition A Hermitian manifold $(M,I,g)$ is
LCK if and only if there exists a closed 1-form $\theta$ (called {\bf
	the Lee form}) such that the fundamental form
$\omega(x,y):=g(Ix,y)$ satisfies
$d\omega=\theta\wedge\omega$.

\hfill

\definition An LCK structure $(\omega,\theta)$ on a complex manifold $(M,I)$ is called {\bf strict} if the Lee form is not exact. When the Lee form is exact, the structure $(\omega,\theta)$ is called {\bf globally conformally K\"ahler} (GCK).

\hfill

\example\label{_Classical_Hopf_Example_} Let  $H:=(\C^n\setminus 0)/\langle A\rangle$, $A=\lambda\Id$, 
$|\lambda|> 1$, be a {\bf classical Hopf manifold}. It is
LCK because the flat K\"ahler metric $\tilde g_0=\sum
dz_i\otimes d\bar z_i$ is multiplied by 2 by the deck
group $\Z$. The Lee form on the covering $\C^n\setminus 0$ 
is $\theta=-d\log |z|^2$. Since the Hopf manifold is diffeomorphic with $S^1\times S^{2n-1}$, the LCK structure is strict. 

\hfill

\remark All complex submanifolds of an LCK manifold are still LCK (but not necessarily strict).

\subsection{LCK manifolds with potential}

\definition\label{_LCK_pote_def_} An LCK manifold has  {\bf LCK potential} if it
admits a K\"ahler covering on which the K\"ahler metric
has a global and positive  potential function $\psi$ such that
the deck group multiplies $\psi$ by a
constant.\footnote{In the sequel, a
	differential form which is multiplied by a constant
	factor by the action of the deck group is called {\bf
		automorphic}.} 
In this case, $M$ is called {\bf LCK manifold with
	potential}.

\hfill

\example The diagonal Hopf manifold in \ref{_Classical_Hopf_Example_} is LCK with potential because $\tilde g_0$ has the global automorphic potential $\psi:=\sum |z_i|^2$. 

By \cite{ov_lckpot}, compact LCK manifolds with potential
are stable to small deformations of the complex
structure. This implies that
all {\bf linear Hopf manifolds} $(\C^n\setminus 0)/\langle A\rangle$, $A\in\mathrm{GL}(n,\C)$,
with eigenvalues of absolute value $> 1$, are LCK with potential.

\hfill

\remark One can easily prove that all complex submanifolds of an LCK manifold with potential are still LCK with potential.

\hfill

One of the most important properties of LCK manifolds with potential is the following Kodaira-type embedding result:

\hfill

\theorem (\cite{ov_lckpot}) Let $M$ be a compact LCK manifold with potential of complex dimension at least 3. Then $M$ admits a holomorphic embedding into a linear Hopf manifold. \endproof

\hfill

The idea of the proof is that  the K\"ahler cover of an LCK
manifold with potential is not a cone in the Riemannian sense
(it lacks the homothety), but is very close to the idea of
cone, because it can be completed with just one (singular)
point. This completion can be realized using a theorem by
Rossi and Andreotti-Siu which only works in dimension
greater than 3. In dimension 2, the result is true,
if one assumes the GSS conjecture (\cite[Theorem 5.6]{ov_indam}), 
which is a long-standing conjecture, often assumed to be true.

Below, we give a precise description of the K\"ahler
cover of an LCK manifold with potential.

\subsubsection{Algebraic cones}

\definition\label{_alge_cone_Definition_} 
(\cite{ov_pams}) 
A {\bf closed algebraic cone}  is an
affine variety $\cac$ admitting a $\C^*$-action $\rho$
with a unique fixed point $x_0$, called {\bf the origin} or {\bf the apex},
and satisfying the following: 
\begin{description}
\item[(i)] $\cac$ is smooth outside of  $x_0$, and $\C^*$ acts freely 
on $\cac \backslash x_0$.
\item[(ii)] $\rho$ acts on the Zariski tangent
space $T_{x_0}\cac$ with all eigenvalues
$|\alpha_i|<1$.
\end{description}
An {\bf open algebraic cone} is a complex manifold
which is biholomorphic to a closed algebraic
cone without the origin.

\hfill

\theorem (\cite{ov_ac})\label{_cone_cover_for_LCK_pot_Theorem_}
Let $M$ be an LCK manifold with proper potential,
and $\tilde M$ its K\"ahler $\Z$-cover.
Then $\tilde M$ is an open algebraic cone. \endproof

\hfill

\remark The restriction to proper potentials is not significant, since 
any LCK metric with potential can be approximated in the ${C}^\infty$-topology  (while keeping the complex structure fixed) by an LCK metric with proper potential (\cite{{ov_jgp_09}, ov_jgp_16}).

\hfill

As a topological space, the closed algebraic cone
$\tilde M_c$ is obtained as a one-point completion of an open
algebraic cone $\tilde M$. However, the complex structure
on a closed algebraic cone is not uniquely determined by
$\tilde M$.  We refer the reader to \cite{ov_ac} for this
intricate and fascinating subject.
 

The manifold $\tilde M$ does not define
the LCK manifold $M$ uniquely: $M= \frac{\tilde M}{\Z}$ 
is determined by the choice of the $\Z$-action which
may vary. 


In \cite{ov_ac}, we proved that an algebraic cone of an
LCK manifold with potential
is isomorphic to an affine cone over a projective orbifold.
This nearly transforms the geometry of LCK manifolds with
potential into a chapter of algebraic geometry.

\hfill

We end these preliminaries by citing a result which gives a method of constructing LCK manifolds with potential once the algebraic cone which will be the covering is known. We first need the following notion.

\hfill 

\definition 
Let $\tilde M$ be an open algebraic cone, $\tilde M_c$
the corresponding closed cone, 
and $\vec r\in T\tilde M_c$ a holomorphic vector field
such that for all $t>0$ the diffeomorphism $e^{t\vec r}$
is a holomorphic contraction of $\tilde M_c$ to the origin.
In this situation, the 1-parameter 
family $e^{t\vec r}$ is called {\bf the flow of contractions} (or {\bf contraction flow}).
A strictly pseudoconvex hypersurface $S\subset \tilde M$ is called 
a {\bf  pseudoconvex shell} if $S$ intersects each orbit
of $e^{t\vec r}$, $t\in \R$, exactly once.

\hfill

Pseudoconvex shells are used to obtain 
automorphic plurisubharmonic functions, as follows.

\hfill

\theorem  (\cite{ov_pams})\label{shell_char}
Let $\tilde M$
be an algebraic cone, $e^{t\vec r}$ a contraction flow, and
$S\subset \tilde M$ a pseudoconvex shell. Then for each $\lambda\in \R$
there exists a unique function $\phi_\lambda$ such that 
$\Lie_{\vec r}\phi = \lambda \phi$ and $\phi_\lambda\restrict S=1$.
Moreover,
{such $\phi_\lambda$  is strictly 
	plurisubharmonic when $\lambda \gg 0$.} \endproof

\subsubsection{Vaisman manifolds}
\label{_Vaisman_Subsubsection_}

Among the LCK manifolds with potential, the best understood subclass is formed by the {\bf Vaisman manifolds}. These are LCK manifolds 
with the Lee form  parallel with respect to the Levi-Civita
connection of the LCK metric. 

Let $(M,I,\omega,\theta)$ be a Vaisman manifold and
$\pi\ :\ \tilde M\to M$ a K\"ahler cover on which
$\pi^*\theta$ is exact. Then it can be seen that
the $\tilde\omega$-squared norm of $\pi^*\theta$ is a
positive, automorphic potential for the K\"ahler metric
$\tilde\omega$ (\cite{ov_lckpot}).

\hfill

\example\label{_Vaisman_Examples_}	
Some examples of Vaisman manifolds:
\begin{description}
	\item[(i)] Diagonal Hopf manifolds $(\C^n\backslash 0)/\langle A\rangle$ where $A$ is semi-simple and with eigenvalues $\alpha_i$ of absolute value $>1$, \cite{go, ov_pams}. 
	\item[(ii)] Elliptic complex surfaces (see \cite{bel} for the complete classification of Vaisman compact surfaces; see also \cite{_ovv:surf_}). 
	\item[(iii)] All compact submanifolds of a Vaisman manifold (\cite{_Verbitsky:Vanishing_LCHK_}).
\end{description}

\remark The class of Vaisman manifolds is strict:  neither
the LCK Inoue surfaces, nor the non-diagonal
Hopf 
manifolds can bear Vaisman metrics (\cite{bel}, \cite{ov_pams}). 

\hfill

The metric dual $\theta^\sharp$ of the Lee form is called
{\bf the Lee field}.

\hfill

\proposition\label{_Lee_field_Killing_Corollary_} (\cite{va_rendiconti, va_gd}) Let $(M,g,I,\theta)$ be a  Vaisman manifold. Then the Lee field
$\theta^\sharp$ is Killing and holomorphic; moreover, 
it commutes with $I\theta^\sharp$. \endproof

\hfill

Denote by $\Sigma$ the holomorphic 1-dimensional foliation generated by $\theta^\sharp$ and $I\theta^\sharp$. It is called {\bf the canonical foliation} (the motivation is given in the next theorem).

\hfill

\theorem \label{_Subva_Vaisman_Theorem_}
Let $M$ be a compact Vaisman manifold, and 
$\Sigma\subset TM$ its canonical foliation. Then:
\begin{description}
	\item[(i)] $\Sigma$ is
	independent from the choice of the Vaisman metric (\cite{tsu}).
	\item[(ii)] $d^c\theta=\omega-\theta \wedge I\theta$ (\cite{va_gd}) and the exact (1,1)-form $\omega_0:= d^c\theta$ is semi-positive
	(\cite{_Verbitsky:Vanishing_LCHK_}). Therefore, $\Sigma=\ker\omega_0$,
	and $\omega_0$ is transversally K\"ahler with respect to $\Sigma$.
\end{description}

The following result is almost obvious. We give its proof for consistency.

\hfill

\lemma \label{_invariance_for_cohomology_Vaisman_Lemma_}
Let $M$ be a compact Vaisman manifold, and $G$
the closure of the Lie group generated by the
action of the Lee and the anti-Lee fields on $M$.
Then $G$ is a compact torus.

\hfill

\proof Since the Lee and anti-Lee fields are acting on $M$ 
by isometries (\ref{_Lee_field_Killing_Corollary_}), 
the group $G$ is a closed subgroup of the
isometry group of $M$, and the latter is compact.
Then $G$ is compact; since it is the closure of an
abelian group, it is abelian.
\endproof


\section{Holomorphic tensors on LCK manifolds with potential}


The purpose of this section is to prove
\ref{_LCK_pot_Zariski_closure_Proposition_}.
For its proof, we shall recall some facts related
to coherent sheaves.

\subsection{Extension of reflexive coherent sheaves}

Given a coherent sheaf $\cal F$ over $X$, let
${\cal F}^*:= \Hom(\cal F, \calo_X)$ be {\bf the dual sheaf},
that is, the sheaf of module homomorphisms to the ring
of regular functions. The natural morphism of sheaves
${\cal F}\arrow {\cal F}^{**}$ does not need to be an isomorphism:
for example, its kernel contains the torsion of ${\cal F}$.
A coherent sheaf is called {\bf reflexive}
if the natural map ${\cal F}\arrow {\cal F}^{**}$ is an isomorphism.
For an introduction to the reflexive sheaves, 
see \cite{_oss_}. Here we state some results which
are relevant to our work.

 First of all, notice that
the sheaf ${\cal F}^*$ is already reflexive
(\cite{_oss_}, Ch. II, Lemma 1.1.8).
Moreover, the singular set of a reflexive sheaf
over a normal variety has codimension $\geq 3$;
in particular, a reflexive sheaf over a complex
surface is locally free (\cite{_oss_}, Ch. II, 1.1.10).
Also, a reflexive sheaf of rank one over a smooth manifold
is locally free (\cite{_oss_},  Ch. II, Lemma 1.1.15).

The most important property of reflexive sheaves is {\bf normality}.
Let $Z\subset X$ be a subvariety of a normal 
complex variety, $\codim_X Z\geq 2$.
Consider the open embedding map
$j:\; X \backslash Z \hookrightarrow X$, and let $j^*$ and $j_*$ be the
sheaf pullback and pushforward, For any sheaf ${\cal F}$, there exists
a natural sheaf morphism  ${\cal F}\arrow j_* j^* {\cal F}$ taking a section
of ${\cal F}$ to its restriction to $X \backslash Z$. A coherent sheaf
${\cal F}$ is {\bf normal} if 
the natural map ${\cal F}\arrow j_* j^* {\cal F}$ is an isomorphism, for any 
subvariety $Z\subset X$ of codimension $\geq 2$.

\hfill

The following theorem can be understood as 
a sheaf version of the Hartogs extension
theorem. 

\hfill

\theorem\label{_reflexife_normal_Theorem_}
A coherent sheaf ${\cal F}$ over a normal complex variety
is normal if and only if it is reflexive.

\hfill

\proof
\cite[Proposition 7]{_Serre:Prolongement_}.
\endproof

\hfill

\theorem \label{_extension_over_a_point_Theorem_} 
(\cite[Theorem 27.6]{_OV:book_}) 
Let $M$ be a normal complex variety with isolated singularities,
$\dim_\C M \geq 3$, $x\in M$ a point, 
and $M_0:= M \backslash \{x\}\stackrel j\hookrightarrow M$.
Consider a reflexive sheaf ${\cal F}$ on $M_0$.
Then $j_*{\cal F}$ is a reflexive coherent sheaf.
Moreover, $j_*{\cal F}$ is the unique coherent
reflexive sheaf on $M$ which is isomorphic to ${\cal F}$
on $M_0$.

\hfill

\proof
Uniqueness of a reflexive extension
follows from the fact that a coherent sheaf ${\cal F}$ over a normal variety 
is normal if and only if it is reflexive (\ref{_reflexife_normal_Theorem_}).
The existence of coherent extension of ${\cal F}$ follows from
\cite[Proposition 6.1]{andreotti_siu}.
\endproof

\subsection{Invariance of holomorphic tensor fields}

We are now ready to prove the main result of this section.

\hfill

\proposition\label{_LCK_pot_Zariski_closure_Proposition_}
Let $M$ be a compact LCK manifold with potential,
$\tilde M$ its K\"ahler $\Z$-cover, considered as
an open algebraic cone, and 
$\Phi\in H^0(M, B)$  be a holomorphic tensor field on $M$,
where  $B=(\Omega^1 M)^{\otimes k}\otimes TM^{\otimes l}$. 
Denote by  $\tilde \Phi$
its lift to $\tilde M$, and let ${\cal G}$ be the Zariski closure
of the $\Z$-action on $\tilde M$. Then $\tilde\Phi$ is
${\cal G}$-invariant.

\hfill

\pstep
Consider a $\Z$-action on a finite-dimensional space $W$.
Then any $\Z$-invariant vector $w\in W$ is invariant under the
Zariski closure of $\Z$. This is why the statement
of \ref{_LCK_pot_Zariski_closure_Proposition_}
is not surprising. However, the space of 
tensor fields on $\tilde M$ is not
finite-dimensional, which makes \ref{_LCK_pot_Zariski_closure_Proposition_}
non-trivial. Denote by ${\goth m}$ the maximal ideal
of the origin in the closed algebraic cone $\tilde M_c$.
We think of the quotients $\calo_{\tilde M_c}/ {\goth m}^k$
as of the spaces of jets of holomorphic (or algebraic)
functions. As in \cite{ov_ac},
we choose $\tilde M_c$ normal, so that the $\Z$-action
is extended to $\tilde M_c$. Then any
$\Z$-invariant jet $u \in \calo_{\tilde M_c}/ {\goth m}^k$
is also ${\cal G}$-invariant.

\hfill

{\bf Step 2:} From now on, we use the same letter $\Phi$
to denote the lift of $\Phi$ to $\tilde M$.
Consider $\Phi$ as a
$\Z$-invariant section of the appropriate
tensor bundle $B=(\Omega^1 M)^{\otimes k}\otimes TM^{\otimes l}$,
which is by construction $\Z$-equivariant.
Using \ref{_extension_over_a_point_Theorem_}, we extend
$B$ to a reflexive coherent sheaf $B_c$ on $\tilde M_c$.
Since coherent reflexive sheaves are normal 
(\cite[Ch. II, Lemma 1.1.12]{_oss_}, \cite[Proposition 7]{_Serre:Prolongement_}), the section $\Phi$
admits a holomorphic extension 
to $\tilde M_c$, denoted as  $\Phi_c\in H^0(\tilde M_c, B_c)$.
Consider the space $J_B^k:=\frac{B_c}{{\goth m}^{k+1} B_c}$
of $k$-jets of $B_c$. The $\Z$-action on $H^0(\tilde M_c, B_c)$
preserves 
$H^0(\tilde M_c, {\goth m}^k B_c)\subset H^0(\tilde M_c, B_c)$.
Let $\Phi^k_c \in J_B^k$ be the $k$-jet of $\tilde\Phi_c$.
Since $\Phi_c$ is $\Z$-invariant, and the space of $k$-jets
is finite-dimensional, the $k$-jet $\Phi^k_c$
is ${\cal G}$-invariant. Let 
$\hat B_c:= \lim\limits_{\leftarrow} J_B^k$
be the ${\goth m}$-adic completion of $B_c$, and
$\hat \Phi_c$ the image of $\Phi_c$ in this
completion. Since all $k$-jets $\Phi^k_c$ are ${\cal G}$-invariant,
the completion $\hat \Phi_c$ is also ${\cal G}$-invariant.
Since the ring of germs of holomorphic
functions on $\tilde M_c$ in $c$ is Noetherian 
(\cite[Chapter II, Theorem B.9]{_Gunning_Rossi_}), and $B_c$, 
being reflexive, is torsion-free, the
completion map $H^0(\tilde M_c, B_c)\arrow \hat B_c$
is injective (\cite[Theorem 10.17]{_Atiyah_MacDonald_}). 
Therefore, $\Phi\in  H^0(\tilde M_c, B_c)$
is also ${\cal G}$-invariant.
\endproof


\section{Zariski closures and the Chevalley theorem}


We recall the following famous theorem due to Chevalley.

\hfill

Let $G\subset \GL(V)$ be an algebraic
group, and $W= V^{\otimes k} \otimes (V^*)^{\otimes l}$
a tensor representation of $G$.  A $G$-invariant
vector $v \in W$ is called {\bf a tensor invariant of $G$}.
A point $x\in {\Bbb P} W$ is called
{\bf a projective tensor invariant of $G$}
if it is $G$-invariant. 

\hfill

\theorem
{\bf (Chevalley theorem,} \cite{_Morris:Ratner_}).\\
An algebraic
group is uniquely determined by the set
of its projective tensor invariants.
\endproof

\hfill

When an algebraic group is {\em reductive} (over $\C$, a group
is reductive if and only if it has a compact real form),
a stronger version of  Chevalley theorem is available.

\hfill

\theorem\label{_Chevalley_reductive_Theorem_}
A reductive algebraic
group is uniquely determined by the set
of its tensor invariants.

\proof 
\cite[Proposition 3.1 (c)]{_Deligne:Hodge_cycles_}.
\endproof

\hfill

We are going to prove the following result, used further
on to determine the Lee field action on a Vaisman manifold.

\hfill

\theorem\label{_Zariski_contains_real_part_Theorem_}
Let $A \in \GL(n, \C)$ be a diagonal linear operator
with the eigenvalues $\alpha_1, ..., \alpha_n$,
and $A_1$ an operator which is diagonal in the
same basis, with the eigenvalues $|\alpha_1|,..., |\alpha_n|$.
Denote by ${\cal G}$ the Zariski closure of the
group $\langle A \rangle$ in $\GL(n, \C)$.
Then ${\cal G}$ contains $A_1$.

\hfill

\proof
Since the operator $A$ is diagonalizable, 
the algebraic closure of $\langle A \rangle$ 
is a commutative group which contains only
semisimple elements. By 
\cite[Proposition 1.5]{_Borel_Tits:Groupes_Reductifs_},
the connected component of ${\cal G}$ is isomorphic
to $(\C^*)^d$, hence ${\cal G}$ is reductive.
Then \ref{_Chevalley_reductive_Theorem_}
implies that ${\cal G}$ is the set of all
$g\in \GL(n, \C)$ which preserve all
$A$-invariant vectors 
$w\in W= V^{\otimes k} \otimes (V^*)^{\otimes l}$.
The eigenvalues of $A$ on $W$
are products of $k$ instances of $\alpha_i$
and $l$ instances of $\alpha_i^{-1}$.
The eigenvectors $w \in W$ 
are $
w = \bigotimes_{i=1}^k z_{m_i} \otimes \bigotimes_{i=1}^l \zeta_{n_i}$, where
$z_{n_i}$ are the eigenvectors in $V$ with the
eigenvalues $\alpha_{n_i}$
and $\zeta_{m_i}$ are the eigenvectors in $V^*$ with the
eigenvalues $\alpha_{m_i}^{-1}$. Then
$A(w) = \prod_{i=1}^k \alpha_{n_i} \prod_{i=1}^l \alpha^{-1}_{m_i}$.
The eigenvector $w$ is $A$-invariant
if and only if $\prod _{i=1}^k \alpha_{n_i} \prod_{i=1}^l \alpha^{-1}_{m_i}=1$.
This implies $\prod _{i=1}^k |\alpha_{n_i}| \prod_{i=1}^l |\alpha^{-1}_{m_i}|=1$,
hence $w$ is $A_1$-invariant as well.
\endproof


\section{Holomorphic tensors on Vaisman manifolds}


Let $M$ be a Vaisman manifold, and $\theta^\sharp$ its Lee 
field. Then $\theta^\sharp$ and $I(\theta^\sharp)$
are holomorphic Killing fields (\ref{_Lee_field_Killing_Corollary_}).
The closure $G$ of the group 
$e^{t_1 \theta^\sharp+t_2 I(\theta^\sharp)}$ in the group of isometries of 
$M$ is compact and  commutative, hence it is a compact torus
(\ref{_invariance_for_cohomology_Vaisman_Lemma_}).

\hfill

\theorem\label{_holo_tensor_on_Vaisman_Lee_invariant_Theorem_}
Let $M$ be a compact Vaisman manifold, 
$\Phi\in H^0(M, B)$ 
a holomorphic tensor, where $B= (\Omega^1 M)^{\otimes k}\otimes TM^{\otimes l}$,
and $G$ the smallest closed Lie group containing the
flows generated by the Lee and the anti-Lee fields.
Then $\Phi$ is $G$-invariant. 

\hfill

\pstep
Choose a Vaisman metric on $M$ with LCK rank 1, and
let $\tilde M$ be the corresponding 
K\"ahler $\Z$-cover of $M$, which is considered
as an open algebraic cone. Let ${\cal G}\subset \Aut(\tilde M)$ be 
the Zariski closure of the $\Z$-action. Then $\Phi$ is
${\cal G}$-invariant by \ref{_LCK_pot_Zariski_closure_Proposition_}.
To finish the proof of \ref{_holo_tensor_on_Vaisman_Lee_invariant_Theorem_},
it would suffice to show that $G \subset {\cal G}$;
since $ {\cal G}$ is closed, this would follow if
we prove that the Lee field $\theta^\sharp$
is tangent to ${\cal G}$.

\hfill

{\bf Step 2:} 
Let $M \arrow H$ be a holomorphic embedding of $M$ to a
diagonal Hopf manifold. By \ref{_Subva_Vaisman_Theorem_}, this embedding
commutes with the group generated by the
Lee and the anti-Lee flow, hence it would suffice
to show that $\theta^\sharp$ is tangent to ${\cal G}$ when $M=H$.
It remains to prove  \ref{_holo_tensor_on_Vaisman_Lee_invariant_Theorem_}
assuming that $M$ is a Hopf manifold. For this purpose,
we compute the groups $G$ and ${\cal G}$ for 
a diagonal Hopf manifold. This can be done explicitly
using the expression for the Vaisman metric on a Hopf manifold
obtained in \ref{shell_char} via the pseudoconvex shells.

\hfill

{\bf Step 3:} The rank 2 algebra generated by
the Lee and anti-Lee flows on a Vaisman
manifold $M$ is uniquely determined by 
the complex structure on $M$ (\ref{_Subva_Vaisman_Theorem_}). 
We express the Lee flow on a Hopf manifold
using a Vaisman structure
we shall construct explicitly, and
prove that it belongs to $\Lie({\cal G})$.
This would imply that the original Lee flow 
also belongs to $\Lie({\cal G})$.

Let $H = \frac{\C^n \backslash 0}{{\langle A \rangle}}$,
where $A \in \GL(n, \C)$ is a diagonal linear contraction with eigenvalues
$\alpha_1, ..., \alpha_n$.  Let $A_1$ be the matrix which is
diagonal in the same basis as $A$, and has eigenvalues
$|\alpha_i|$. Then $\vec r:=\log A_1$
is a holomorphic contraction vector field,
and  $e^{t\vec r}$ is an $A$-invariant
contraction flow. The sphere $S^{2n-1}$ 
is a pseudoconvex shell compatible with this contraction flow.
By \cite{ov_pams}, we obtain a $e^{t\log A_1}$-automorphic
plurisubharmonic function $\phi$. We are going to prove
that $\phi$ is $A$-automorphic.

Consider the diagonal $\U(1)^n$-action on $\C^n$.
The plurisubharmonic function $\phi$ is $\U(1)^n$-invariant
because $\vec r$ and $S^{2n-1}$ are $\U(1)^n$-invariant.
Since $A\in \U(1)^n \cdot A_1$,
and $\phi$ is $A_1$-automorphic, it is $A$-automorphic.
Therefore, $\frac{dd^c\phi}{\phi}$ defines a Vaisman metric 
on $H$. 
By definition, the Lee flow maps the level sets
of $\phi$ to the level sets of $\phi$, with 
$\phi(e^{t\theta^\sharp}z)= t\phi(z)$.
By construction of $\phi$,
the same is true for $e^{-t\log A_1}$, that is,
$\phi(e^{-t\log A_1}z)= t\phi(z)$.
Therefore, $\theta^\sharp = -\log A_1$, where $\log A_1$
is the diagonal matrix with the eigenvalues $\log |\alpha_i|$
in the same basis. This is the Lee flow for the
Vaisman metric $\frac{dd^c\phi}{\phi}$. 

The vector field
$\log A_1$ belongs to $\Lie({\cal G})$,
because any power of $A_1$ fixes all
tensor invariants of ${\cal G}$
by \ref{_Zariski_contains_real_part_Theorem_},
hence the corresponding 1-parametric
sugroup also belongs to $\Lie({\cal G})$;
this gives $\log A_1\in \Lie({\cal G})$.
\endproof


\section{An application: the Kodaira dimension of Vaisman manifolds}


\definition
Let $M$ be a compact complex manifold.
Its {\bf pluricanonical bundle} is
the tensor power $K_M^n= K_M^{\otimes n}$, $n \geq 0$.
The {\bf  Kodaira dimension} $\kappa(M)$  
is defined as $\kappa(M):=\limsup_n \frac {\log(\dim H^0(K_M^n))} {\log n}$.

\hfill

\remark
The Kodaira dimension is equal to $-\infty$ if $H^0(K_M^n)=0$
for almost all $n>0$, and to 0 if $\dim H^0(K_M^n)$ is bounded, but non-zero
for infinitely many $n$. 
If the function $n \mapsto \dim H^0(K_M^n)$
grows as a polynomial of degree $d$, the Kodaira dimension
of $M$ is $d$.

\hfill

Let $M$ be a Vaisman manifold, 
$G$ the closure of the group generated by the
Lee and the anti-Lee action, $\tilde M$ 
a K\"ahler $\Z$-covering of $M$, and $\tilde G$ the lift
of $G$ to $\tilde M$. Denote by ${\goth G}$ the Zariski
closure of $\tilde G$ in $\Aut(\tilde M)$, considered
as an algebraic variety (\cite{ov_ac}). 

Since $\tilde G$ acts on $\tilde M$ by homotheties,
and contains contractions, there is an open
set $U\subset {\goth G}$ in the connected component
of ${\goth G}$ such that for all  $\gamma'\in U$,
the quotient $M' := \frac{\tilde M}{\langle \gamma'\rangle}$
is LCK (\cite{ov_ac}). The automorphism $\gamma$ commutes with the action
of the group generated by the Lee and the anti-Lee
field, hence by \cite{kor} $M'$ is Vaisman.\footnote{As
  shown in  \cite{kor}, an LCK manifold which admits
a holomorphic conformal $\C$-action is Vaisman, assuming
that this action cannot be lifted to an isometry of the
K\"ahler cover.}
By \cite[Proposition 1.5]{_Borel_Tits:Groupes_Reductifs_},
the connected component of ${\goth G}$ is
isomorphic to $(\C^*)^k$.
It is not hard to see that the union of all one-dimensional
closed subgroups is dense in any algebraic group, hence
$U$ contains a dense subset of $\gamma'$ which are
contained in a subgroup $\C^* \subset {\goth G}$.
Since $\C^*/\langle \gamma' \rangle$ is an elliptic
curve, the Vaisman manifold 
$M' = \frac{\tilde M}{\langle \gamma'\rangle}$
is elliptically fibered.\footnote{Ellipticaly fibered Vaisman manifolds
are called {\bf quasi-regular}.} As shown in
\cite{ov_imm_vai}, the leaf space $X$ is a 
projective orbifold.

\hfill

\theorem\label{_Kodaira_dimension_Theorem_}
Let $M$ be a Vaisman manifold,
$M'$ its quasi-regular deformation obtained
above, and $X$ the corresponding leaf space,
which is a projective orbifold. Then
$\kappa(M)=\kappa(M')=\kappa(X)$,
where $\kappa$ denotes the Kodaira dimension.

\hfill

\pstep The equality
$\kappa(M')=\kappa(X)$ is clear from the adjunction
formula: since $\pi:\; M' \arrow X$ is an elliptic
fibration, one has 
$H^0(K_{M'}^{k})=H^0(K_{X}^{ k})$ for all $k$.

\hfill

{\bf Step 2:}
By \ref{_holo_tensor_on_Vaisman_Lee_invariant_Theorem_},
any section $h$ of the pluricanonical bundle $K_M^k$ on $M$ 
is ${\goth G}$-invariant. Then, for any
$\gamma'\in U$, the lift of $h$ to $\tilde M$
is $\gamma'$-invariant, hence can be obtained
as a pullback of a section of $K_{M'}^k$ under
the quotient map $\tilde M \arrow M'$.
This defines an injective map
$\Psi:\; H^0(K_M^k) \arrow H^0(K_{M'}^k)$.

\hfill

{\bf Step 3:}  To finish the proof
of \ref{_Kodaira_dimension_Theorem_} it remains
to show that $\Psi:\; H^0(K_M^k) \arrow H^0(K_{M'}^k)$
is surjective. This would follow if we show that
any section $h' \in H^0(K_{M'}^k)$
is ${\goth G}$-invariant. Since 
$H^0(K_{M'}^{k})=H^0(K_{X}^{k})$ 
(Step 1), the ${\goth G}$-invariance
of $h'$ would follow if we prove that
the corresponding section of $K_{X}^{k}$ is
invariant under the natural action
of ${\goth G}$ on $X$. 

Recall that {\bf the pluricanonical representation}
is the natural action of $\Aut(X)$ on $H^0(K_{X}^{k})$.
Since $X$ is an orbifold, its singularities are klt
(\cite[Proposition 1.2.1]{_Prokhorov:log_surface_}).
By \cite[Theorem 10.61]{_Kollar_Kovacz_},
the image of the pluricanonical representation
is finite; since ${\goth G}$ is connected,
it acts on $H^0(K_{X}^{k})$ trivially.
This gives $H^0(K_{M}^{k})=H^0(K_{M'}^{k})=H^0(K_{X}^{ k})$ .
\endproof

\hfill

As an application, we  outline
the proof of the deformational stability
of Kodaira dimension for Vaisman manifolds.

\hfill

\conjecture\label{_Kodaira_constant_Conjecture_}
Let $M_t, t \in \R$, be a smooth family of compact
Vaisman manifolds. Then the Kodaira dimension
of $M_t$ is constant.

\hfill

{\bf Sketch of a proof:}
We sketch a proof of this conjecture, reducing it
to \ref{_cone_smoothly_Conjecture_} below.

Let $\tilde M_t$ be the family of open algebraic
cones obtained from each $M_t$ by taking
a K\"ahler $\Z$-cover; this is possible to do by
\cite{ov_imm_vai}. Each $\tilde M_t$ is equipped
with a natural algebraic structure by \cite{ov_ac}.

\hfill

What follows is a conjecture, which most likely 
can be proven by the methods used in \cite{ov_ac}.

\hfill

Let $\tilde M_t, t \in \R$,
be a smooth family of open algebraic cones. 
As shown in \cite{ov_ac}, for an
appropriate choice of the $\C^*$-action on $\tilde M_t$
the quotient $\tilde M_t/\C^*$ is isomorphic to a
complex projective orbifold $X_t$; this isomorphism
is extended to an isomorphism between $\tilde M_t$ 
and $\Tot^\circ(L_t)$,
where $L_t$ is an ample line  bundle over $X_t$.

\hfill

\conjecture \label{_cone_smoothly_Conjecture_}
Let $\tilde M_t, t \in \R$,
be a smooth family of open algebraic cones.
Then for any sufficiently small open
subset $U\subset \R$, the choice of $\C^*$-action 
as above can be performed in such a way that
 the family $(X_t, L_t)$ 
depends smoothly on $t\in U$. 

\hfill

Now, by \ref{_Kodaira_dimension_Theorem_},
the Kodaira dimension of $M_t$ is equal
to the Kodaira dimension of $X_t$, which
is constant, as shown by Siu in
\cite{_Siu:plurigenera_}.
\endproof

\hfill

{\bf Acknowledgements:}
We are grateful to Nicolina Istrati for her interest and 
valuable advice.

\hfill

\hfill

{\small

\noindent {\sc Liviu Ornea\\
{\sc University of Bucharest, Faculty of Mathematics and Informatics, \\14
Academiei str., 70109 Bucharest, Romania}, and:\\
 Institute of Mathematics ``Simion Stoilow" of the Romanian
Academy,\\
21, Calea Grivitei Str.
010702-Bucharest, Romania}\\
{\tt lornea@fmi.unibuc.ro,   liviu.ornea@imar.ro}

\hfill

\noindent
{\sc Misha Verbitsky\\
{\sc Instituto Nacional de Matem\'atica Pura e
	Aplicada (IMPA) \\ Estrada Dona Castorina, 110\\
	Jardim Bot\^anico, CEP 22460-320\\
	Rio de Janeiro, RJ - Brasil }\\
also:\\
Laboratory of Algebraic Geometry, \\
Faculty of Mathematics, National Research University 
Higher School of Economics,
6 Usacheva Str. Moscow, Russia}\\
\tt verbit@verbit.ru, verbit@impa.br

 }

\end{document}